\documentclass[10pt]{amsart}
\input{diagrams}

\newtheorem{proposition}{Proposition}[section]
\newtheorem{lemma}[proposition]{Lemma}

\newtheorem{theorem}[proposition]{Theorem}

\theoremstyle{definition}
\newtheorem{definition}[proposition]{Definition}
\newtheorem{example}[proposition]{Example}

\theoremstyle{remark}
\newtheorem{remark}[proposition]{Remark}

\newtheorem{blanco}[proposition]{}

\newcommand{\thlabel}[1]{\label{th:#1}}
\newcommand{\thref}[1]{Theorem~\ref{th:#1}}
\newcommand{\selabel}[1]{\label{se:#1}}
\newcommand{\seref}[1]{Section~\ref{se:#1}}
\newcommand{\lelabel}[1]{\label{le:#1}}
\newcommand{\leref}[1]{Lemma~\ref{le:#1}}
\newcommand{\prlabel}[1]{\label{pr:#1}}
\newcommand{\prref}[1]{Proposition~\ref{pr:#1}}

\newcommand{\relabel}[1]{\label{re:#1}}
\newcommand{\reref}[1]{Remark~\ref{re:#1}}
\newcommand{\exlabel}[1]{\label{ex:#1}}
\newcommand{\exref}[1]{Example~\ref{ex:#1}}
\newcommand{\delabel}[1]{\label{de:#1}}
\newcommand{\deref}[1]{Definition~\ref{de:#1}}
\newcommand{\eqlabel}[1]{\label{eq:#1}}
\newcommand{\equref}[1]{(\ref{eq:#1})}

\newcommand{\Hom}{{\rm Hom}}

\newcommand{\Ker}{{\rm Ker}\,}
\newcommand{\Coker}{{\rm Coker}\,}

\def\ot{\otimes}

\def\mapright#1{\smash{\mathop{\longrightarrow}\limits^{#1}}}
\def\mapleft#1{\smash{\mathop{\longleftarrow}\limits^{#1}}}

\newcommand{\Cc}{\mathcal{C}}
\newcommand{\Dd}{\mathcal{D}}

\newcommand{\Ff}{\mathcal{F}}

\newcommand{\Mm}{\mathcal{M}}

\newcommand{\Tt}{\mathcal{T}}

\newcommand{\Ww}{\mathcal{W}}
\newcommand{\Zz}{\mathcal{Z}}
\newcommand{\Ss}{\mathcal{S}}
\newcommand{\Xx}{\mathcal{X}}
\newcommand{\YD}{\mathcal{YD}}
\def\*C{{}^*\hspace*{-1pt}{\Cc}}
\def\text#1{{\rm {\rm #1}}}

\def\ol{\overline}
\def\ul{\underline}
\def\dul#1{\underline{\underline{#1}}}

\def\doublerightbis#1#2{{\lower.2ex\vbox{
\hbox{${\smash{\mathop{\longrightarrow}\limits^{#1}}}$}\vspace*{-4mm}
\hbox{${\smash{\mathop{\longrightarrow}\limits_{#2}}}$}}}}

\begin{document}
\title[Turaev's
Hopf group-coalgebras]{A categorical approach to Turaev's
Hopf group-coalgebras}
\author{S. Caenepeel}
\address{Faculty of Engineering Sciences,
Vrije Universiteit Brussel, VUB, B-1050 Brussels, Belgium}
\email{scaenepe@vub.ac.be}
\urladdr{http://homepages.vub.ac.be/\~{}scaenepe/}
\author{M. De Lombaerde}
\address{Faculty of Engineering Sciences,
Vrije Universiteit Brussel, VUB, B-1050 Brussels, Belgium}
\email{mdelomba@vub.ac.be}
\thanks{}


\maketitle

\begin{abstract}
We show that Turaev's group-coalgebras and Hopf group-coalgebras are
coalgebras and Hopf algebras in a symmetric monoidal category, which we
call the Turaev category. A similar result holds for group-algebras
and Hopf group-algebras. As an application, we give an alternative
approach to Virelizier's version of the Fundamental Theorem for Hopf
algebras. We introduce Yetter-Drinfeld modules over Hopf group-coalgebras
using the center construction.
\end{abstract}

\section*{Introduction}
Group-coalgebras and Hopf group-coalgebras appeared in the work of Turaev
\cite{Turaev} on homotopy quantum field theories. In the case where
the underlying group $G$ is trivial, we recover the classical coalgebras
and Hopf algebras. A purely algebraic study of Hopf group-coalgebras
was initiated by Virelizier \cite{Virelizier}, and then continued by
Zunino \cite{Zunino1,Zunino2} and Wang \cite{Wang1,Wang2,Wang3}.
It turns out that many of the classical results in Hopf algebra theory
can be generalized to the group coalgebra setting. Virelizier
gives a generalized version of the Fundamental Theorem for Hopf algebras,
and introduces $G$-integrals; Zunino introduces Yetter-Drinfeld modules,
the Drinfeld double, and a generalization of the center construction of
a monoidal category; Wang introduces Doi-Hopf modules, entwined modules
and coalgebra Galois theory for Hopf group-coalgebras, and he proves
a version of Maschke's Theorem.\\
In contrast with classical Hopf algebras, the definition of a Hopf group-coalgebra
is not selfdual. In fact, there exist the dual notions of  group-algebra
and Hopf group-algebra. The dual of a group-coalgebra is a group-algebra,
and the converse property holds under some finiteness assumptions.\\
In this paper, we propose an alternative approach: we will introduce
a symmetric monoidal category $\Tt_k$, called the Turaev category, and show
that coalgebras (resp. Hopf algebras) in $\Tt_k$ are precisely
group-coalgebras (resp. Hopf group-coalgebras). The objects of $\Tt_k$
are $k$-vector spaces (or $k$-modules, if $k$ is a commutative ring),
indexed by a set $X$. The morphisms are defined in such a way that we have
a strongly monoidal forgetful functor to the opposite of the category of
sets. We can define a second monoidal category $\Zz_k$, called the Zunino
category, with the same objects, but differently defined morphisms,
such that the algebras (resp. Hopf algebras) in $\Zz_k$ are the
group-algebras (resp. Hopf group-algebras). The forgetful functor is now
a strongly monoidal functor to the category of sets, and this will explain
the lack of selfduality in the definition of Hopf group-coalgebras.\\
Some of the results of the above cited papers can now be viewed from the
perspective of symmetric monoidal categories. We will discuss two
examples. In \seref{3}, we will show that Virelizier's Fundamental Theorem
can be viewed as a special case of Takeuchi's version of the Fundamental
Theorem in braided monoidal categories \cite{Takeuchi}. In \seref{4},
we will compute the center of the category of modules over a Hopf group-coalgebra,
which will lead to the introduction of Yetter-Drinfeld modules.

\section{Preliminary results}\selabel{1}
\subsection{Monoidal categories}\selabel{1.1}
Recall from e.g. \cite{Kassel} or \cite{McLane} that a monoidal category (or a
tensor category) consists of a triple $\Cc=(\Cc,\ot, I)$, where
$\Cc$ is a category, $\ot:\ \Cc\times \Cc\to \Cc$ is a bifunctor, and $I$ is an object
in $\Cc$, together with natural isomorphisms
$$(U\ot V)\ot W\cong U\ot (V\ot W)$$
$$V\ot I\cong V\cong I\ot V$$
for all $U,V,W\in \Cc$, satisfying some appropriate coherence conditions.
In all our examples, the associativity constraint is the obvious one, so we will identify
$(U\ot V)\ot W$ and $U\ot (V\ot W)$, and we will just write $U\ot V\ot W$.
In a similar way, we identify $I\ot V$, $V$ and $V\ot I$.\\
A braided monoidal category is a monoidal category together with a natural 
isomorphism $\tau_{V,W}:\ V\ot W\to W\ot V$ called the braiding or commutativity constraint,
satisfying some appropriate coherence conditions. The category is called symmetric
if $\tau_{W,V}$ is the inverse of $\tau_{V,W}$, for all $V,W$ in $\Cc$.\\
An algebra (or monoid) in a monoidal category $\Cc$ is a triple $(A,m,\eta)$,
with $A$ an object in $\Cc$, and $m:\ A\ot A\to A$ and $\eta:\ I\to A$
morphisms in $\Cc$ such that
$$m\circ (A\ot m)= m\circ (m\ot A),~~m\circ (\eta\ot A)=m\circ(A\ot \eta)= A,$$
where $A$ denotes the identity morphism of $A$. A right $A$-module is
a couple $(M,\psi)$, with $M\in \Cc$ and $\psi:\ M\ot A\to M$ such that
$$\psi\circ (\psi\ot A)=\psi\circ(M\ot m)~~{\rm and}~~\psi\circ M\ot \eta=M$$
The category of right $A$-modules will be denoted by $\Mm_A$.\\
The tensor product of two algebras $A$ and $B$ in a braided monoidal category
is again an algebra, with multiplication
$$m_{A\ot B}=(m_A\ot m_B)\circ (A\ot \tau_{B,A}\ot B).$$
If $\Cc$ is a (braided) monoidal category, then $\Cc^{\rm op}$ is also a (braided) monoidal category,
with the same tensor product and unit, and with the original constraints replaced by
their inverses. A coalgebra in $\Cc$ is an algebra in $\Cc^{\rm op}$.\\
If $(C,\Delta,\varepsilon)$ is a coalgebra, and $A$ is an algebra, then $\Hom_{\Cc}(C,A)$ is a semigroup,
with multiplication (the convolution) given by
$$f*g=m\circ (f\ot g)\circ \Delta$$
The unit for the convolution is $\eta\circ \varepsilon$.\\
A bialgebra in a braided monoidal category
$\Cc$ is a fivetuple $H=(H,m,\eta,\Delta,\varepsilon)$ such that
$(H,m,\eta,\Delta)$ is an algebra and $(H,\Delta,\varepsilon)$ is a coalgebra,
and such that $\eta$ and $m$ are are coalgebra maps, or, equivalently,
$\delta$ and $\varepsilon$ are algebra maps. If the identity on $H$ has an
inverse $S$ in $\Hom_\Cc(H,H)$, then $H=(H,m,\eta,\Delta,\varepsilon, S)$
is called a Hopf algebra in $\Cc$. $S$ is called the antipode.

\begin{example}\exlabel{1.3}
Let $k$ be a commutative ring, and consider the category $\Mm_k$ of $k$-modules
and $k$-linear maps. Then $(\Mm_k,\ot_k, k)$ is a symmetric monoidal category.
A algebra in this category is a usual $k$-algebra. The same is true for
a coalgebra, a bialgebra and a Hopf algebra. A module over an algebra $A$ is
a usual $A$-module, and the same is true for a comodule over a coalgebra.
\end{example}

\begin{example}\exlabel{1.4}
Let $\dul{\rm Sets}$ be the category of sets and functions; $X \times Y$
is the cartesian product of two sets $X$ and $Y$; fix a singleton $\{*\}$. Then
$(\dul{\rm Sets},\times,\{*\})$ is a symmetric monoidal category. An algebra in $\dul{\rm Sets}$
is a semigroup.
\end{example}

The following result is folklore, but we give details, as it will be useful
in the sequel.

\begin{lemma}\lelabel{1.5}
Let $X$ be a set. There exists a unique comultiplication $\Delta$ and counit 
$\varepsilon$ on $X$ making $(X,\Delta,\varepsilon)$ a coalgebra in 
$\dul{\rm Sets}$. $\Delta$ and $\varepsilon$ are the following maps:
\begin{equation}\eqlabel{1.5.1}
\Delta:\ X\to X\times X,~~\Delta(x)=(x,x)~~;~~
\varepsilon:\ X\to\{*\},~~\varepsilon(x)=*.
\end{equation}
A bialgebra $(X,m,\eta,\Delta,\varepsilon)$ in 
$\dul{\rm Sets}$ is of the following form: $(X,m,\eta)$ is a semigroup,
and $\Delta$ and $\varepsilon$ are given by \equref{1.5.1}.\\
Finally, the category of Hopf algebras in $\dul{\rm Sets}$ is isomorphic
to the category of groups.
\end{lemma}

\begin{proof}
Let $(X,\Delta,\varepsilon)$ be a coalgebra in $\dul{\rm Sets}$.
Then $\varepsilon$ is the unique map from $X$ to $\{*\}$. Take
$x\in X$, and assume that $\Delta(x)=(y,z)\in X\times X$. Then
$$(x,*)=(y,\varepsilon(z))~~{\rm and}~~(*,x)=(\varepsilon(y),z)$$
hence $x=y=z$. Conversely, it is a straightforward to check that
$(X,\Delta,\varepsilon)$ is a coalgebra in $\dul{\rm Sets}$.
If $(X,m,\eta)$ is a semigroup, then the maps $\Delta$ and
$\varepsilon$ are semigroup maps, hence $(X,m,\eta,\Delta,\varepsilon)$
is a bialgebra in $\dul{\rm Sets}$.\\
Let $X$ be a set (hence a coalgebra), and $G$ a semi-group (hence an
algebra). The convolution on ${\rm Map}(X,G)$ is then just pointwise
multiplation: for $f,g:\ X\to G$, we have
$$(f* g)(x)=f(x)g(x)$$
The unit element of the convolution is $\eta\circ \varepsilon:\ X\to G$,
mapping every $x\in X$ to the unit element $1$ of $G$.\\
Now let $(X,m,\eta,\Delta,\varepsilon, s)$ be a Hopf algebra in $\dul{\rm Sets}$.
Then $s*X=X*s= \eta\circ \varepsilon$, hence for all $x\in X$:
$$(m\circ (s,X) \circ \Delta)(x)= (m\circ (X,s) \circ \Delta)(x)=1$$
or
$$s(x)x=xs(x)=1$$
so $s(x)=x^{-1}$, and $X$ is a group. Conversely, if $X$ is a group, then we
can give it a bialgebra structure, and $s(x)=x^{-1}$ is an antipode.
\end{proof}

Let $\Cc$ and $\Dd$ be monoidal categories. A monoidal functor 
$F:\ \Cc\to \Dd$ is
a functor together with a natural transformation
$$u_{V,W}:\ F(V)\ot F(W)\to F(V\ot W)$$
and a morphism $$J\to F(I)$$
compatible with the constraints. If this natural transformation and this
morphism are isomorphisms, then we call the monoidal functor strong.
A monoidal functor sends algebras to algebras, and a strong monoidal functor
also sends coalgebras to coalgebras. A braided monoidal functor between two {braided}
mocategories is a monoidal functor between the monoidal categories that is
compatible with the braidings, i.e.
$$u_{W,V}\circ t_{F(V),F(W)}=F(t_{V,W})\circ u_{V,W}$$
for all $V,W\in \Cc$. A braided strong monoidal functor sends bialgebras to
bialgebras and Hopf algebras to Hopf algebras.\\

A strong monoidal functor sends algebras, coalgebras, bialgebras
and Hopf algebras to algebras, coalgebras, bialgebras and Hopf algebras.
An example of a strong monoidal functor is
the functor $F:\ \dul{\rm Sets}\to \Mm_k$ sending a set $X$ to the free $k$-module
with basis $X$. $F$ sends a group $G$ (a Hopf algebra in $\dul{\rm Sets}$) to the
group algebra $kG$, which is a Hopf algebra in $\Mm_k$.\\

If $A$ is an algebra in $\Cc$, then a right $A$-module $M$ is an object $M\in \Cc$
together with a morphism $\psi:\ M\ot A\to M$ such that
$$\psi\circ (\psi\ot A)=\psi\circ (A\ot \mu)~~{\rm and}~~
\psi\circ (M\ot \eta)=M$$
A morphism $f$ between two right $A$-modules $M$ and $N$ is called right $A$-linear
if $\psi\circ (f\ot A)=f\circ \psi$. The category of right $A$-modules is denoted
by $\Cc_A$.\\
The category of comodules $\Cc^C$ over a coalgebra $C$ is defined in a dual fashion.
If $H$ is a bialgebra in a braided monoidal category $\Cc$, a (right) Hopf module is
an object $M\in \Cc$ together with $\psi:\ M\ot H\to M$ and $\rho:\ M\to M\ot H$
such that
\begin{equation}\eqlabel{1.11}
\rho\circ\psi= (\psi\ot \mu)\circ t_{23}\circ (\rho\ot \Delta)
\end{equation}

\subsection{The center of a monoidal category}\selabel{1.1a}
Let $\Cc$ be a monoidal category. The weak left center $\Ww_l(\Cc)$ is the category
with the following objects and morphisms. An object is a couple $(V,s)$,
with $V\in \Cc$ and $s:\ V\ot -\to -\ot V$  a natural transformation between
functors $\Cc\to \Cc$, satisfying the following condition, for all
$X,Y\in \Cc$:
\begin{equation}\eqlabel{center1}
(X\ot s_Y)\circ (s_X\ot Y)=s_{X\ot Y},
\end{equation}
and such that $s_I$ is the composition of the natural isomorphisms
$V\ot I\cong V\cong I\ot V$. A morphism between $(V,s)$ and $(V',s')$
consists of $f:\ V\to V'$ in $\Cc$ such that
$$(X\ot f)\circ s_X=s'_X\circ (f\ot X)$$
$\Ww_l(\Cc)$ is a prebraided monoidal category. The unit is $(I,I)$, and
the tensor product is
$$(V,s)\ot (V',s')=(V\ot V', u)$$
with
\begin{equation}\eqlabel{center2}
u_X=s_X\ot V'\circ V\ot s'_X
\end{equation}
The prebraiding $c$ on $\Ww_l(\Cc)$ is given by
\begin{equation}\eqlabel{center3}
c_{V,V'}=s_{V'}:\ (V,s)\ot (V',s')\to (V',s')\ot (V,s).
\end{equation}
The left center $\Zz_l(\Cc)$ is the full subcategory of $\Ww_l(\Cc)$
consisting of objects $(V,s)$ with $s$ a natural isomorphism. $\Zz_l(\Cc)$
is a braided monoidal category.
For detail in the case where $\Cc$ is a strict monoidal category, we refer
to \cite[Theorem XIII.4.2]{Kassel}. The results remain valid in the case of
an arbitrary monoidal category, since every monoidal category is equivalent
to a strict one.

\subsection{$ G $-coalgebras and Hopf $ G $-coalgebras}\selabel{1.2}
From \cite{Turaev}, we recall the following definitions.

\begin{definition}\delabel{1.6}
Let $ G $ be a group, en $k$ a commutative ring. A $ G $-coalgebra
is a family of $k$-modules $C=(C_ g )_{ g \in  G }$ indexed by
the group $ G $, together with a family of linear maps
$$\Delta=(\Delta_{ g , h }:\ C_{ g  h }\to C_{ g }\ot
C_{ h })_{ g , h \in  G }$$
and a linear map $\varepsilon:\ C_1\to k$ such that the following conditions
hold, for all $ g , h , k \in  G $:
\begin{equation}\eqlabel{1.6.1}
(\Delta_{ g , h }\ot C_ k )\circ \Delta_{ g  h , k }=
(C_ g \ot \Delta_{ h , k })\circ \Delta_{ g , h  k },
\end{equation}
\begin{equation}\eqlabel{1.6.2}
(C_ g \ot\varepsilon)\circ \Delta_{ g ,1}=C_ g =
(\varepsilon\ot C_ g )\circ \Delta_{1, g },
\end{equation}
\end{definition}

In \cite{Virelizier}, the following Sweedler-type notation for comultiplication
is introduced: for $c\in C_{gh}$, we write
$$\Delta_{g,h}(c)=c_{(1,g)}\ot c_{(2,h)}.$$
For $c\in C_{ghk}$, we will write, using the coassociativity,
$$((\Delta_{ g , h }\ot C_ k )\circ \Delta_{ g  h , k })(c)=
((C_ g \ot \Delta_{ h , k })\circ \Delta_{ g , h  k })(c)=
c_{(1,g)}\ot c_{(2,h)}\ot c_{(3,k)}.$$

\begin{definition}\delabel{1.7}
We use the notation from \deref{1.6}. A semi-Hopf $ G $-coalgebra
is a family of $k$-algebras $H=(H_ g )_{ g \in  G }$
indexed by the group $ G $, together with families of linear maps $\Delta$
and $\varepsilon$ such that $H$ is a $ G $-coalgebra, and such that
$\varepsilon$ and $\Delta_{ g , h }$ are algebra maps.\\
A Hopf $ G $-coalgebra is a semi-Hopf $ G $-coalgebra together with a
family of maps
$$S=\{S_ g :\ H_ g ^{-1}\to H_{ g }~|~ g \in  G \}$$
such that
\begin{equation}\eqlabel{1.7.1}
\mu_ g \circ(S_{ g }\ot H_ g )\circ \Delta_{ g ^{-1}, g }=
\eta_ g \circ \varepsilon= \mu_{ g }\circ (H_ g \ot S_{ g })
\circ \Delta_{ g ^{-1}, g },
\end{equation}
for all $ g \in  G $, where $\mu_ g $ and $\eta_ g $ are the multiplication
and counit maps on $H_ g $.
\end{definition}

Observe that a $ G $-coalgebra and a semi-Hopf $ G $-coalgebra
can also be defined in the case where
$ G $ is just a monoid. 

\subsection{$ G $-algebras and Hopf $ G $-algebras}\selabel{1.3}
The definition of $ G $-coalgebra and Hopf $ G $-coalgebra can be dualized.
This was already remarked in \cite{Turaev}; a formal definition has been
presented in \cite{Zunino1}.

\begin{definition}\delabel{1.8}
Let $ G $ be a group. A $ G $-algebra consists of a set of $k$-modules
$A=(A)_{ g \in G }$ together with maps
$$\mu_{ g , h }:\ A_ g \ot A_ h \to A_{ g  h }~~{\rm and}~~
\eta:\ k\to A_1$$
such that
\begin{equation}\eqlabel{1.8.1}
\mu_{ g  h , k }\circ (\mu_{ g , h }\ot A_ k )=
\mu_{ g , h  k }\circ (A_ g \ot \mu_{ h ,k})
\end{equation}
and
\begin{equation}\eqlabel{1.8.2}
\mu_{ g ,1}\circ (A_ g \ot \eta)=A_ g =
\mu_{1, g }\circ(\eta\ot A_ g )
\end{equation}
\end{definition}

\begin{definition}\delabel{1.9}
Let $ G $ be a group. A $ G $-algebra $H=(H)_{ g \in G }$ is called
a semi-Hopf $ G $-algebra if every $H_ g $ is a $k$-coalgebra in such a
way that $\mu_{ g , h }$ and $\eta$ are $k$-coalgebra maps.
\end{definition}

\begin{definition}\delabel{1.10}
Let $ G $ be a group. A Hopf $ G $-algebra is a semi-Hopf 
$ G $-algebra $H=(H)_{ g \in G }$ together with maps $S_ g :\
H_ g \to H_{ g ^{-1}}$ such that
\begin{equation}\eqlabel{1.10.1}
\mu_{ g ^{-1}, g }\circ (S_ g \ot H_ g )\circ\Delta_ g =
\mu_{ g , g ^{-1}}\circ (H_ g \ot S_ g )\circ\Delta_ g =
\eta\circ\varepsilon_ g 
\end{equation}
\end{definition}

Note that $pi$-algebras and semi-Hopf $ G $-algebras can be defined in
the situation where $ G $ is only a semigroup.

\section{The Turaev and Zunino categories}\selabel{2}
In \deref{1.6}, we fix a group $ G $. The central idea is to replace
$ G $ by a variable set.

\subsection{The Turaev category}\selabel{2.1}
\begin{definition}\delabel{2.1}
Let $k$ be a commutative ring. A Turaev $k$-module is
a couple $\ul{M}=(X,M)$, where $X$ is a set, and $M=(M_x)_{x\in X}$ is
a family of $k$-modules indexed by $X$. A morphism between
two $T$-modules $(X,M)$ and $(Y,N)$ is a couple $\ul{\varphi}=(f,\varphi)$,
where $f:\ Y\to X$ is a function, and 
$\varphi=(\varphi_y:\ M_{f(y)}\to N_y)_{y\in Y}$ is a family
of linear maps indexed by $Y$. The composition of $\ul{\varphi}:\ \ul{M}\to \ul{N}$
and $\ul{\psi}:\ \ul{N}\to \ul{P}=(Z,P)$ is defined as follows:
$$\ul{\psi}\circ\ul{\varphi}=(f\circ g, (\psi_z\circ\varphi_{g(z)})_{z\in Z})$$
The category of Turaev $k$-modules is
denoted by $\Tt_k$.
\end{definition}
We will use the following notation for the composition of morphisms:
\begin{equation}\eqlabel{2.1.1}
\begin{matrix}
\ul{M}&\mapright{\ul{\varphi}}&\ul{N}&\mapright{\ul{\psi}}&\ul{P}\\
X&\mapleft{f}&Y&\mapleft{g}&Z\\
M_{f(g(z))}&\mapright{\varphi_{g(z)}}&N_{g(z)}&\mapright{\psi_z}&P_z
\end{matrix}
\end{equation}
The category of $T$-modules is a symmetric monoidal category. The tensor
product of $(X,M)$ and $(Y,N)$ is given by
$$(X,M)\ot (Y,M)=(X\times Y, (M_x\ot N_y)_{(x,y)\in X\times Y}),$$
and the unit object is $(\{*\},k)$. The symmetry
 $$\ul{\tau}=(t,\tau):\
(X\times Y, (M_x\ot N_y)_{(x,y)\in X\times Y})\to
(Y\times X, (N_y\ot M_x)_{(y,x)\in Y\times X})$$
is defined as follows: $t:\ Y\times X\to X\times Y$ is the switch map,
and $\tau_{t(y,x)}:\ M_x\ot N_y\to N_y\ot M_X$ is also the switch map.
We have two strong monoidal functors
$$F':\ \Mm_k\to \Tt_k,~~F'(M)=(\{*\},M)$$
and
$$F:\ \Tt_k\to \dul{\rm Sets}^{\rm op},~~F(X, M)=X$$

\begin{proposition}\prlabel{2.2}
If $( G ,C)$ is a coalgebra in $\Tt_k$,
then $ G $ is a monoid and
$C$ is a $ G $-coalgebra. Conversely, if $C$ is a $ G $-coalgebra,
then $( G ,C)$ is a coalgebra in $\Tt_k$.
\end{proposition}

\begin{proof}
Let $\ul{C}=(\ul{C}=( G ,C),\ul{\Delta}=(m,\Delta),\ul{\varepsilon}=(i,\varepsilon))$ be a coalgebra in $\Tt_k$.
Then $F(\ul{C})=( G ,m,i)$ is a coalgebra in $\dul{\rm Sets}^{\rm op}$, hence a
monoid. The counit and comultiplication are
\begin{equation}\eqlabel{2.2.1}
\begin{matrix}
\ul{C}&\mapright{\ul{\varepsilon}}&\ul{k}\\
 G &\mapleft{i}&\{*\}\\
C_1=C_{i(e)}&\mapright{\varepsilon}&k
\end{matrix}~~~~{\rm and}~~~~
\begin{matrix}
\ul{C}&\mapright{\ul{\Delta}}&\ul{C}\ot\ul{C}\\
 G &\mapleft{m}& G \times G \\
C_{ g  h }=C_{m( g , h )}&\mapright{\Delta_{ g , h }}&C_{ g }
\ot C_{ h }
\end{matrix}
\end{equation}
We have
$$\begin{matrix}
\ul{C}&\mapright{\ul{\Delta}}&\ul{C}\ot\ul{C}&\mapright{\ul{C}\ot\ul{\Delta}}&
\ul{C}\ot\ul{C}\ot\ul{C}\\
 G &\mapleft{m}& G \times G &\mapleft{( G ,m)}& G \times G \times G \\
C_{ g  h  k }&\mapright{\Delta_{ g , h  k }}&C_{ g }
\ot C_{ h  k }& \mapright{C_ g \ot \Delta_{ h , k }}&
C_{ g }\ot C_{ h }\ot C_{ k }
\end{matrix}$$
and
$$\begin{matrix}
\ul{C}&\mapright{\ul{\Delta}}&\ul{C}\ot\ul{C}&\mapright{\ul{\Delta}\ot\ul{C}}&
\ul{C}\ot\ul{C}\ot\ul{C}\\
 G &\mapleft{m}& G \times G &\mapleft{(m, G )}& G \times G \times G \\
C_{ g  h  k }&\mapright{\Delta_{ g  h , k }}&C_{ g  h }
\ot C_{ k }& \mapright{\Delta_{ g , h }\ot C_ k }&
C_{ g }\ot C_{ h }\ot C_{ k }
\end{matrix}$$
hence $\ul{\Delta}$ is coassociative if and only if \equref{1.6.1} holds,
for all $ g , h , k \in  G $. In a similar way
$$\begin{matrix}
\ul{C}&\mapright{\ul{\Delta}}&\ul{C}\ot\ul{C}&\mapright{\ul{\varepsilon}
\ot\ul{C}}&\ul{C}\\
 G &\mapleft{m}& G \times G &\mapleft{(i, G )}& G \\
C_{ g }&\mapright{\Delta_{1, g }}&C_{1}
\ot C_{ g }& \mapright{\varepsilon\ot C_ g }&C_ g 
\end{matrix}$$
and consequently $({\ul{\varepsilon}\ot\ul{C}})\circ \ul{\Delta}=
\ul{C}$ if and only if the second equality in \equref{1.6.2} holds,
for all $ g \in  G $. In a similar way, 
$({\ul{C}\ot \ul{\varepsilon}})\circ \ul{\Delta}=
\ul{C}$ if and only if the first equality in \equref{1.6.2} holds.
Hence $C$ is a $ G $-coalgebra.\\
Conversely, given a $ G $-coalgebra $C$, we define a coalgebra structure
on $( G ,C)$ by \equref{2.2.1}.
\end{proof}

\begin{proposition}\prlabel{2.3}
An algebra $(X,A)$ in $\Tt_k$ consists of a set $X$ and a family
of $k$-algebras $A=(A_x)_{x\in X}$ indexed by $X$. A map $(f,\varphi)$
in $\Tt_k$
between two algebras 
$(X,A)$ and $(Y,B)$ is an algebra map if and only if every
$\varphi_{f(y)}:\ A_{f(y)}\to B_y$ is a $k$-algebra map.
\end{proposition}

\begin{proof}
Take an algebra $\ul{A}=((X,A),\ul{\mu}=(\delta,\mu),\ul{\eta}=
(e,\eta))$ in $\Tt_k$. Then $F(\ul{A})$ is an algebra in 
$\dul{\rm Sets}^{\rm op}$, hence a coalgebra in $\dul{\rm Sets}$.
It follows from \leref{1.5} that $X$ is an arbitray set, and
that $\delta(x)=(x,x)$ and $e(x)=*$. The multiplication and unit are
morphisms
$$\begin{matrix}
\ul{A}\ot \ul{A}&\mapright{\ul{\mu}}&\ul{A}\\
X\times X&\mapleft{\delta}&X\\
A_x\ot A_x&\mapright{\mu_x}& A_x
\end{matrix}
~~~{\rm and}~~~
\begin{matrix}
\ul{k}&\mapright{\ul{\eta}}&\ul{A}\\
\{*\}&\mapleft{e}&X\\
k&\mapright{\eta_x}& A_x
\end{matrix}$$
Let us compute
$$\begin{matrix}
\ul{A}\ot \ul{A}\ot \ul{A}&\mapright{\ul{A}\ot\ul{\mu}}&
\ul{A}\ot \ul{A}&\mapright{\ul{\mu}}&\ul{A}\\
X\times X\times X&\mapleft{(X,\delta)}&
X\times X&\mapleft{\delta}&X\\
A_x\ot A_x\ot A_x&\mapright{A_x\ot\mu_x}&
A_x\ot A_x&\mapright{\mu_x}& A_x
\end{matrix}$$
In a similar way, we compute $\ul{\mu}\circ (\ul{A}\ot \ul{\mu})$, and it
follows that $\ul{\mu}$ is associative if and only if every $\mu_x$ is
associative.\\
We also have
$$\begin{matrix}
\ul{A}&\mapright{\ul{A}\ot\ul{\eta}}&\ul{A}\ot\ul{A}&\mapright{\ul{\mu}}&
\ul{A}\\
X&\mapleft{(X,e)}&X\times X& \mapleft{\delta}& X\\
A_x&\mapright{A_x\ot \eta_x}&A_x\ot A_x&\mapright{\mu_x}&A_x
\end{matrix}$$
hence $\ul{\mu}\circ(\ul{A}\ot \ul{\eta})=\ul{A}$ if and only if
$\mu_x\circ (A_x\ot \eta_x)=A_x$ for every $x\in X$. The same is true
for the left unit property, and our result follows. The final statement is
straightforward, and is left to the reader.
\end{proof}

\begin{proposition}\prlabel{2.4}
If $( G ,H)$ is a bialgebra in $\Tt_k$,
then $ G $ is a monoid and
$H$ is a semi-Hopf $ G $-coalgebra. Conversely, if $H$ is a semi-Hopf $ G $-coalgebra,
then $( G ,H)$ is a bialgebra in $\Tt_k$.
\end{proposition}

\begin{proof}
Let $( G ,H)$ be a bialgebra in $\Tt_k$. We know from Propositions
\ref{pr:2.2} and \ref{pr:2.3} that $ G $ is a monoid, $H$ is a $ G $-coalgebra,
and every $H_ g $ is a $k$-algebra. It follows from the final statement in
\prref{2.4} that $\ul{\varepsilon}$ and $\ul{\Delta}$ are algebra maps if
and only if $\varepsilon$ and $\Delta_{ g , h }$ are $k$-algebra maps,
for all $ g , h \in  G $.
\end{proof}

\begin{proposition}\prlabel{2.5}
If $( G ,H)$ is a Hopf algebra in $\Tt_k$,
then $ G $ is a group and
$H$ is a Hopf $ G $-coalgebra. Conversely, if $H$ is a Hopf $ G $-coalgebra,
then $( G ,H)$ is a Hopf algebra in $\Tt_k$.
\end{proposition}

\begin{proof}
Let $\ul{H}=( G ,H)$ be a Hopf algebra in $\Tt_k$. Then $F(\ul{H})= G $ is
a Hopf algebra in $\dul{\rm Sets}^{\rm op}$, and also in $\dul{\rm Sets}$, since
the definition of a Hopf algebra in a category is selfdual. Thus $ G $ is a group.
Let $s:\  G \to  G $, $s( g )= g ^{-1}$, and consider a map
$\ul{S}=(s,S):\ \ul{H}\to \ul{H}$ in $\Tt_k$. We compute the convolution
$\ul{S}*\ul{H}$.
$$\begin{matrix}
\ul{H}&\mapright{\ul\Delta}&\ul{H}\ot \ul{H}&\mapright{\ul{S}\ot\ul{H}}&
\ul{H}\ot \ul{H}&\mapright{\ul{\mu}}&\ul{H}\\
 G &\mapleft{m}& G \times G &\mapleft{(s, G )}& G \times G & \mapleft{\delta}& G \\
H_1&\mapright{\Delta_{ g ^{-1}, g }}&H_{ g ^{-1}}\ot H_{ g }&
\mapright{S_{ g }\ot H_{ g }}& H_{ g }\ot H_{ g }&\mapright{\mu_ g }&
H_ g 
\end{matrix}$$
We also compute 
$$\begin{matrix}
\ul{H}&\mapright{\ul{\varepsilon}}&\ul{k}& \mapright{\ul{\eta}}&\ul{H}\\
 G &\mapleft{1}&\{*\}& \mapleft{e} &G \\
H_1&\mapright{\varepsilon}&k&\mapright{\eta_ g }&H_ g 
\end{matrix}$$
So $\ul{S}*\ul{H}= \ul{\eta}\circ\ul{\varepsilon}$ if and only if
$$\mu_ g \circ(S_{ g }\ot H_ g )\circ \Delta_{ g ^{-1}, g }=
\eta_ g \circ \varepsilon$$
In a similar way, $\ul{H}*\ul{S}= \ul{\eta}\circ\ul{\varepsilon}$ if and oly if
$$\mu_{ g }\circ (H_ g \ot S_{ g })
\circ \Delta_{ g ^{-1}, g }=\eta_ g \circ \varepsilon,$$
and our result follows, in view of \equref{1.7.1}
\end{proof}

The next result will turn out to be important when we discuss the Fundamental
Theorem for Hopf $ G $-coalgebras. But let us first recall that the category
$\dul{\rm Sets}$ has coequalizers (see e.g. \cite[p. 65]{McLane}). For two maps
$f,g:\ Y\to X$, one considers the smallest equivalence relation $\sim$
on $X$ containing
$$\{(f(y),g(y))~|~y\in Y\}$$
The coequalizer of $f$ and $g$ is the natural surjection $x\to \ol{X}=X/\sim$.
The class in $\ol{X}$ represented by $x\in X$ will be denoted by $\ol{x}$,
and observe that $f^{-1}(\ol{x})=g^{-1}(\ol{x})$.

\begin{proposition}\prlabel{2.6}
The category $\Tt_k$ has equalizers and coequalizers.
\end{proposition}

\begin{proof}
Take two morphisms in $\Tt_k$:
$$\begin{matrix}
\ul{M}&\mapright{\ul{\varphi}}&\ul{N}\\
X&\mapleft{f}&Y\\
M_{f(y)}&\mapright{\varphi_y}&N_y
\end{matrix}~~{\rm and}~~~~~
\begin{matrix}
\ul{M}&\mapright{\ul{\psi}}&\ul{N}\\
X&\mapleft{g}&Y\\
M_{g(y)}&\mapright{\psi_y}&N_y
\end{matrix}$$
Take $\ol{x}\in \ol{X}$, and put
$$\ol{M}_{\ol{x}}=\{(m_x)\in \prod_{x\in\ol{x}}M_x~|~\forall y\in f^{-1}(\ol{x})=g^{-1}(\ol{x}):~~\varphi_y(m_{f(y)})=\psi_y(m_{g(y)})\}$$
Set $\ul{\ol{M}}=(\ol{X},\ol{M}=(\ol{M}_{\ol{x}})_{\ol{x}\in\ol{X}})$. 
We claim that the equalizer of
$\ol{\varphi}$ and $\ol{\psi}$ is
$$\begin{matrix}
\ul{\ol{M}}&\mapright{}& \ul{M}\\
\ol{X}&\mapleft{}&X\\
\ol{M}_{\ol{x}}&\mapright{}&M_x
\end{matrix}$$
where $\ol{M}_{\ol{x}}\to M_x$ is the restriction to $\ol{M}_{\ol{x}}$
of the projection $\prod_{x\in\ol{x}}M_x\to M_x$.\\
Consider a morphism $\ul{\pi}:\ \ul{P}\to\ul{M}$ in $\Tt_k$ such that
$\ul{\varphi}\circ\ul{\pi}=\ul{\psi}\circ\ul{\pi}$. Thus
$$\begin{matrix}
\ul{P}&\mapright{\ul{\pi}}&\ul{M}&\mapright{\ul{\varphi}}&\ul{N}\\
Z&\mapleft{h}&X&\mapleft{f}&Y\\
P_{h(f(y))}&\mapright{\pi_{f(y)}}&M_{f(y)}&\mapright{\varphi_y}&N_y
\end{matrix}~~{\rm and}~~
\begin{matrix}
\ul{P}&\mapright{\ul{\pi}}&\ul{M}&\mapright{\ul{\psi}}&\ul{N}\\
Z&\mapleft{h}&X&\mapleft{g}&Y\\
P_{h(g(y))}&\mapright{\pi_{g(y)}}&M_{g(y)}&\mapright{\psi_y}&N_y
\end{matrix}$$
We now define a morphism $\ol{\ul{\pi}}:\ \ul{P}\to\ul{\ol{M}}$. The map
$$\ol{h}:\ \ol{X}\to Z,~~\ol{h}(\ol{x})=h(x)$$
is well-defined. Take $\ol{x}\in\ol{X}$. For all $x'\in \ol{x}$,
we have that $h(x')=h(x)$, and we have a map $\pi_{x'}:\ P_{x'}\to M_x$.
Now define
$$\begin{matrix}
\ul{P}&\mapright{\ul{\ol{\pi}}}&\ul{\ol{M}}\\
Z&\mapleft{\ol{h}}&\ol{X}\\
P_{\ol{h}(\ol{x})}&\mapright{\ul{\pi}_{\ul{x}}}&\ul{M}_{\ul{x}}
\end{matrix}$$
as follows: we know that $P_{\ol{h}(\ol{x})}=P_{h(x')}$, and we define
$$\ul{\pi}_{\ul{x}}(p)=\Bigl(\pi_{x'}(p)\Bigr)_{x'\in\ul{x}}.$$
Suppose that $x=f(y)$ and $x'=g(y')$. Then 
$$\varphi_y(\pi_{f(y)}(p))=\psi_y(\pi_{g(y)}(p)),$$
hence $\ul{\pi}_{\ul{x}}(p)\in \ol{M}_{\ul{x}}$, as needed.\\
Let us now prove that $\Tt_k$ has coequalizers. Let
$$U=\{u\in Y~|~f(u)=g(u)\}$$
For $u\in U$, we have maps
$$\varphi_u,\psi_u:\ M_{f(u)}=M_{g(u)}\to N_u$$
We now put $P_u=\Coker (\varphi_u-\psi_u)$, en
$\ul{P}=(Z,(P_u)_{u\in U})$.
$$\begin{matrix}
\ul{N}&\mapright{}&\ul{P}\\
Y&\mapleft{}&U\\
M_{f(y)}=M_{g(y)}&\mapright{}&P_u
\end{matrix}$$
is a morphism in $\Zz_k$, and is the coequalizer of $\ul{\varphi}$ and
$\ul{\psi}$.
\end{proof}

\subsection{The Zunino category}\selabel{2.2}
Let $k$ be a commutative ring.
The objects of the Zunino category $\Zz_k$ are the same as the objects of
$\Tt_k$: a set $X$ together with a family of $k$-modules $(M_x)_{x\in X}$ indexed by $X$.
A morphism between $\ul{M}=(X,M= (M_x)_{x\in X})$ and
$\ul{N}=(T,N= (N_y)_{y\in Y})$ consists of a map $f:\ X\to Y$ and a family
of $k$-module homomorphisms $\varphi_x:\ M_x\to N_{f(x)}$. We use the following
notation: $\ul{\varphi}=(f,\varphi):\ \ul{M}\to\ul{N}$, or
$$\begin{matrix}
\ul{M}&\mapright{\ul{\varphi}}&\ul{N}\\
X&\mapright{f}&Y\\
M_x&\mapright{\varphi_x}&N_{f(x)}
\end{matrix}$$
The category $\Zz_k$ is a symmetric monoidal category; the tensor product and unit
object are the same as the ones on $\Tt_k$. We also have strong monoidal functors
$$F':\ \Mm_k\to \Zz_k,~~F'(M)=(\{*\},M)$$
$$F:\ \Zz_k\to \dul{\rm Sets},~~F(X,M)=X$$
The proof of the following results is similar to the corresponding proof of 
Propositions \ref{pr:2.2}, \ref{pr:2.3}, \ref{pr:2.4} and \ref{pr:2.5}. We
leave the details to the reader.

\begin{proposition}\prlabel{2.7}
Assume that $( G ,A)$ is an algebra in $\Zz_k$. Then $ G $ is a monoid, and
$A$ is a $ G $-algebra, in the sense of \deref{1.8}. Conversely, if
$A$ is a $ G $-algebra, then $( G ,A)$ is an algebra in $\Zz_k$.
\end{proposition}

\begin{proposition}\prlabel{2.8}
A coalgebra in $\Zz_k$ consists of a family of coalgebras $C_x$ indexed by
a set $X$.
\end{proposition}

\begin{proposition}\prlabel{2.9}
Assume that $( G ,H)$ is a bialgebra in $\Zz_k$. Then $ G $ is a monoid, and
$H$ is a semi-Hopf $ G $-algebra, in the sense of \deref{1.9}. Conversely, if
$H$ is a semi-Hopf $ G $-algebra, then $( G ,H)$ is a bialgebra in $\Zz_k$.
\end{proposition}

\begin{proposition}\prlabel{2.10}
Assume that $( G ,H)$ is a Hopf algebra in $\Zz_k$. Then $ G $ is a group, and
$H$ is a Hopf $ G $-algebra, in the sense of \deref{1.10}. Conversely, if
$H$ is a Hopf $ G $-algebra, then $( G ,H)$ is a Hopf algebra in $\Zz_k$.
\end{proposition}

Let $k$ be a commutative ring. We have braided monoidal contravariant functors
$${}^*:\ \Tt_k\to \Zz_k~~{\rm and}~~{}^*:\ \Zz_k\to \Tt_k$$
For $\ul{M}=(X,(M_x)_{x\in X}\in \Tt_k$, we let
$$\ul{M}^*=(X,(M^*_x)_{x\in X}.$$
A similar definition holds for $\ul{M}=(X,(M_x)_{x\in X}\in \Zz_k$. Let us
illustrate that the functors are monoidal: for $\ul{M}=\ul{N}\in \Tt_k$,
we have the following morphism in $\Zz_k$:
$$\begin{matrix}
\ul{M}^*\ot \ul{N}^*&\mapright{}&(\ul{M}\ot\ul{N})^*\\
X\times Y&\mapleft{=}& X\times Y\\
M_x^*\ot N_y^*&\mapright{}&(M_x\ot N_y)^*
\end{matrix}$$
Consequently, if $\ul{C}$ is a coalgebra in $\Tt_k$, then $\ul{C}^*$ is
an algebra in $\Zz_k$ (in fact the duality functor is a covariant functor
from $\Tt_k^{\rm op}$ to $\Zz_k$, hence it sends algebras to algebras,
see \seref{1.1}).\\
We call $\ul{M}=(X,(M_x)_{x\in X})\in \Tt_k$ (resp. in $\Zz_k$) {\sl finite}
if $M_x$ is finitely generated and projective as a $k$-module. If
$\ul{M}$ and $\ul{N}$ are finite, then $M_x^*\ot N_y^*\cong (M_x\ot N_y)^*$,
for every $x\in X$ and $y\in Y$. Hence the duality functor induces a
strong monoidal functor between the full subcategories $\Tt_k^f$ and
$\Zz_k^f$ consisting of finite objects, and the two duality functors establish
a pair of inverse equivalences.\\
Thus the dual of a finite algebra in $\Zz_k$ is a coalgebra in $\Tt_k$, and
the dual of a bialgebra or Hopf algebra in $\Tt_k$ (resp. $\Zz_k$) is
a bialgebra or Hopf algebra in $\Zz_k$ (resp. $\Tt_k$). These facts are known
(see \cite{Wang2,Zunino1}).

\begin{proposition}\prlabel{2.11}
The category $\Zz_k$ has equalizers and coequalizers.
\end{proposition}

\begin{proof}
Take two morphisms in $\Zz_k$:
$$\begin{matrix}
\ul{M}&\mapright{\ul{\varphi}}&\ul{N}\\
X&\mapright{f}&Y\\
M_{x}&\mapright{\varphi_x}&N_{f(x)}
\end{matrix}~~{\rm and}~~
\begin{matrix}
\ul{M}&\mapright{\ul{\psi}}&\ul{N}\\
X&\mapright{g}&Y\\
M_{x}&\mapright{\psi_x}&N_{g(y)}
\end{matrix}$$
Let $Z=\Ker(f,g)=\{x\in X~|~f(x)=g(x)\}$, and, for each $x\in Z$,
$P_x=\Ker(\varphi_x-\psi_x)$. Then the equalizer of $\ul{\varphi}$ and
$\ul{\psi}$ is
$$\begin{matrix}
\ul{P}&\mapright{}&\ul{M}\\
Z&\mapright{\subset}&X\\
P_x&\mapright{\subset}&M_x
\end{matrix}$$
Let $\ol{Y}=\Coker(f,g)$. Take $\ol{y}\in \ol{Y}$ and $y'\in\ol{y}$. We have
a map
$$i_{y'}:\ N_{y'}\to \prod_{y\in\ol{y}}N_y,~~
(i_{y'}(n))_y=
\left\{
\begin{array}{ccc}
n&{\rm if}& y=y'\\ 0&{\rm if} &y\neq y'
\end{array}\right.
$$
On $\displaystyle{\prod_{y\in\ol{y}}}N_y$, we consider the smallest equivalence
relation $\sim$ containing
$$\{\Bigl((i_{f(x)}\circ\varphi_x)(m), (i_{g(x)}\circ\psi_x)(m)\Bigr)~|~
x\in f^{-1}(\ol{y})=g^{-1}(\ol{y}),~m\in M_x\},$$
and we define
$$\ol{N}_{\ol{y}}=\prod_{y\in\ol{y}}N_y/\sim.$$
The coequalizer of $\ul{\varphi}$ and
$\ul{\psi}$ is
$$\begin{matrix}
\ul{N}&\mapright{}&\ol{\ul{N}}\\
Y&\mapright{}&\ol{Y}\\
N_y&\mapright{p_y}&\ol{N}_{\ol{y}}
\end{matrix}$$
where $p_y$ is the composition of $i_y$ and the canonical projection
$$\prod_{y\in\ol{y}}N_y\to \prod_{y\in\ol{y}}N_y/\sim.$$
\end{proof}

\section{The Fundamental Theorem}\selabel{3}
It is known, see e.g. \cite[Theorem 3.4]{Takeuchi} that the Fundamental
Theorem of Hopf modules, as stated in Sweedler's book \cite{Sweedler},
can be formulated in any braided monoidal category with equalizers. If the
category has coequalizers, then we have a second version of the Fundamental
Theorem, since a Hopf algebra is also a Hopf algebra in the opposite category.
Since
$\Zz_k$ and $\Tt_k$ have equalizers and coequalizers, we have 
two versions of the Fundamental Theorem in both
categories. The aim of this Section is to make this explicit, and to
derive Virelizier's version of the Fundamental Theorem from it
(see \cite[Theorem 2.7]{Virelizier}). First we look at what is happening in
the category of sets.

\subsection{Hopf modules in the category of sets}\selabel{3.1}
Let $ G $ be a set, and consider the corresponding coalgebra $( G ,\delta,e)$
in $\dul{\rm Sets}$. Let $(X,\rho)$ be a right $ G $-comodule. It is easy to
show that $\rho(x)=(x,f(x))$, with $f:\ X\to G $ an arbitrary map. Thus
a right $ G $-comodule consists of a set $X$ together with a map $f:\ X\to  G $.\\
Now take a semigroup $ G $, that is an algebra in $\dul{\rm Sets}$.
A right $ G $-module is a right $ G $-set $X$, this is a set $X$ together with
a $ G $-action such that $x1=x$ and $(x g ) h =x( g  h )$, for all $x\in X$
and $ g , h \in  G $.\\
A semigroup $ G $ is also a bialgebra in $\dul{\rm Sets}$, so we can consider
right Hopf modules over $ G $. These are right $ G $-sets $X$ together with a
map $f:\ X\to  G $ satisfying
\begin{equation}\eqlabel{3.1.1}
f(x g )=f(x) g 
\end{equation}
for all $x\in X$ and $ g \in G $.\\

Let $X$ be a right $ G $-Hopf module. The coinvariants $X^{{\rm co} G }$ of $X$ are defined as the equalizer of $(X,f)$ and $(X,\eta)$, that is
$$X^{{\rm co} G }=\{x\in X~|~f(x)=1\}$$
The Fundamental Theorem then takes the following form.

\begin{theorem}\thlabel{3.1}
Let $ G $ be a group, and $X$ be a right $ G $-Hopf module. Then the map
$$\phi:\ X^{{\rm co} G }\times  G \to X,~~\phi(x, g )=x g $$
is bijective.
\end{theorem}

\begin{proof}
This is a special case of \cite[Theorem 3.4]{Takeuchi}, but it can also be
verified directly. The inverse of $\phi$ is given by
$$\phi^{-1}(x)=(xf(x)^{-1},f(x))$$
\end{proof}

A group $ G $ is also a Hopf algebra in $\dul{Sets}^{\rm op}$. Let $X$ be a
$ G $-Hopf module, and consider $X_ G $, the coequalizer of
$$X\times  G \mapright{\rm action}~~{\rm and}~~X\times  G  \mapright{(X,e)} X$$
In $X_ G $, $x$ and $x g $ are identified, which means that $X_ G $
consists of the orbits under the $ G $-action. The orbit containing ${x}$
will be denoted by $\ol{x}$. The second version of the Fundamental Theorem
is now the following.

\begin{theorem}\thlabel{3.2}
Let $ G $ be a group, and $X$ be a right $ G $-Hopf module. Then the map
$$\psi:\ X\to X\times  G ,~~\psi(x)=(\ol{x},f(x))$$
is bijective.
\end{theorem}

\begin{proof}
Again a special case of \cite[Theorem 3.4]{Takeuchi}; but we can also verify
directly that $\psi^{-1}$ is given by
$$\psi^{-1}(\ol{x}, g )=xf(x)^{-1} g $$
\end{proof}

\subsection{Modules, comodules and Hopf modules in $\Tt_k$}.\selabel{3.2}
Let $\ul{C}=(G,C)$ be a coalgebra in $\Tt_k$. We can consider the category
$\Tt^{\ul{C}}$ of right $\ul{C}$-comodules and right $\ul{C}$-colinear
maps. An object $\ul{M}=(X,M)\in\Tt^{\ul{C}}$ consists of a right $G$-set $X$
and a set of $k$-modules indexed by $X$, together with a coaction
$$\ul{\rho}:\ \ul{M}\to \ul{M}\ot \ul{C}$$
in $\Tt_k$ such that
$$(\ul{M}\ot\ul{\varepsilon})\circ \ul{\rho}=\ul{M}~~{\rm and}~~
(\ul{\rho}\ot\ul{\Delta})\circ \ul{\rho}=(\ul{\rho}\ot\ul{C})\circ \ul{\rho}$$
Explicitely,
$$\begin{matrix}
\ul{M}&\mapright{\ul{\rho}}& \ul{M}\ot \ul{C}\\
X&\mapleft{}&X\times G\\
M_{xg}&\mapright{\rho_{x,g}}&M_x\ot C_g
\end{matrix}$$
satisfies
$$(M_x\ot\varepsilon)\circ\rho_{x,1}=M_x~~{\rm and}~~
(\rho_{x,g}\ot C_h)\circ \rho_{xg,h}=
(M_x\ot\Delta_{g,h})\circ \rho_{x,gh}$$
Note that our definition is more general than the one in \cite{Virelizier}, where
 $G$ is a group, and $X=G$.

\begin{example}\exlabel{3.2.0}
Let $X$ be a right $G$-set, and $M$ a $k$-module. Let $\rho_g:\ M\to M\ot C_g$
be maps such that
$$(\rho_g\ot C_h)\circ \rho_h=(M\ot \Delta_{g,h})\circ \rho_{gh}~~{\rm and}~~
(M\ot \varepsilon)\circ\rho_1=M.$$
Then $(\{*\},M)$ is a right $\ol{C}$-comodule. For all $x\in X$, let $M_x=M$
as a $k$-module, and $\ul{M}=\{X,(M_x)_{x\in X}\in \Tt_k$. Define
$\ul{\rho}:\ \ul{M}\to \ul{M}\ot\ul{C}$ as follows:
$$\rho_{x,g}=\rho_g:\ M_{xg}=M\to M_x\ot C_g= M\ot C_g.$$
Then $\ul{M}$ is a right $\ul{C}$-comodule.
\end{example}

Now let $\ul{A}=(Y,A)$ be an algebra in $\Tt_k$. A right $\ul{A}$-module
$\ul{M}=(X,M)$ consists of a set of modules $M_x$ indexed by the set $X$,
and a map $f:\ X\to Y$ such that $M_x$ is a right $A_{f(x)}$-module.\\
If $\ul{H}$ is a bialgebra in $\Tt_k$, then we can consider (right)
Hopf modules. Let $\ul{M}=(X,M)$ be such a Hopf module. Then $\ul{M}$ is
a right $\ul{H}$-module and a right $\ul{H}$-comodule, as above. The map
$f:\ X\to G$ has to be right $G$-linear, and \equref{1.11} has to be
satisfied, explicitely
$$(\psi_x\circ M_g)\circ (M_x\ot t_{g,f(x)}\ot M_g)\circ
(\rho_{x,g}\ot \Delta_{f(x),g})=\rho_{x,g}\circ \psi_{xg}$$
for all $x\in X$ and $g\in G$. Here $\psi_x:\ M_x\ot H_{f(x)}\to M_x$
is the right $H_{f(x)}$-action on $M_x$, and $t_{g,f(x)}:\
H_g\ot H_{f(x)}\to H_{f(x)}\ot H_g$ is the switch map.

\subsection{The Fundamental Theorem in the category $\Tt_k$}\selabel{3.3}
Let $\ul{H}$ be a Hopf algebra in $\Tt_k$, and consider a right Hopf module
$\ul{M}$. The module of coinvariants $\ul{M}^{{\rm co}\ul{H}}$ is by definition
the equalizer of $\ul{\rho}$ and $\ul{M}\ot\ul{\eta}$. 
$\ul{M}^{{\rm co}\ul{H}}$ can be described as follows:
$$\ul{M}^{{\rm co}\ul{H}}=(X_G, (M_{\ol x})_{\ol{x}\in X_G})$$
with
$$M^{{\rm co}\ul{H}}_{\ol{x}}=\{(m_y)\in \prod_{y\in \ol{x}}M_y~|~
\rho_{y,g}(m_{xg})=m_x\ot 1_g,~{\rm for ~all~}y\in \ol{x}~{\rm and}~
g\in G\}$$

\begin{theorem}\thlabel{3.3} {\bf (Fundamental Theorem)}
Let $\ul{H}$ be a Hopf algebra in $\Tt_k$, and $\ul{M}$ a Hopf module.
Then we have an isomorphism of Hopf modules
$$\begin{matrix}
\ul{M}^{{\rm co}\ul{H}}\ot \ul{H}&\mapright{\ul{\phi}}& \ul{M}\\
X_G\times G&\mapleft{\varphi}& X\\
M_{\ol{x}}\ot H_{f(x)}&\mapright{\phi_{x,g}}&M_x
\end{matrix}$$
with $\varphi(x)=(\ol{x},f(x))$ and $\phi_{x,g}((m_y)_{y\in\ol{x}})\ot h_{f(x)})=
m_xh_{f(x)}$.
\end{theorem}

\begin{proof}xxx
This is a direct consequence of \thref{3.2}. Let us give the explicit
formula for $\ul{\phi}^{-1}$. First observe that the two compositions
below coincide:
$$\ul{M}\rTo^{\ul{\rho}}\ul{M}\ot\ul{H}\rTo^{\ul{M}\ot\ul{S}}
\ul{M}\ot\ul{H}\rTo^{\ul{\psi}}\ul{M} 
\pile{\rTo^{\ul{\rho}}\\ \rTo_{\ul{M}\ot\ul{H}}}
\ul{M}\ot\ul{H}$$
Thus we have a map $\ul{R}:\ \ul{M}\to \ul{M}^{{\rm co}\ul{H}}$.
According to the proof of \prref{2.6}, $\ul{R}$ can be described explicitely
as follows:
$$\begin{matrix}
\ul{M}&\mapright{\ul{R}}& \ul{M}^{{\rm co}\ul{H}}\\
X&\mapleft{r} & X_G\\
M_{xf(x)^{-1}}&\mapright{R_{\ol{x}}}& M^{{\rm co}\ul{H}}_{\ol{x}}
\end{matrix}$$
with $r(\ol{x})=xf(x)^{-1}$. For $x'\in \ol{x}$, we consider
$m\in M_{x'f(x')^{-1}}=M_{xf(x)^{-1}}$, and
$$\rho_{x',f(x')^{-1}}(m)=m_{[0,x']}\ot m_{[1,f(x')^{-1}]}$$
Then
$$R_{\ol{x}}(m)=\left(m_{[0,x']}S_{f(x')}(m_{[1,f(x')^{-1}]})\right)_{x'\in X}$$
and $\ul{\phi}^{-1}$ is the composition $(\ul{R}\ot\ul{H})\circ \ul{\rho}$,
that is
$$\begin{matrix}
\ul{M}&\mapright{\ul{\phi}^{-1}}&\ul{M}^{{\rm co}\ul{H}}\ot\ul{H}\\
X&\mapleft{\varphi^{-1}}& X_G\times G\\
M_{xf(x)^{-1}g}&\mapright{\phi^{-1}_{\ol{x},g}}&M^{{\rm co}\ul{H}}_{\ol{x}}\ot
H_g
\end{matrix}$$
with
$$\varphi^{-1}(\ol{x},g)=xf(x)^{-1}g$$
and
$$\phi^{-1}_{\ol{x},g}(m)= 
\left(m_{[0,x']}S_{f(x')}(m_{[1,f(x')^{-1}]})\right)_{x'\in X}
\ot m_{[2,g]}$$
\end{proof}

\begin{remark}\relabel{3.4}
Let us consider the particular situation where $X=G$ and $f:\ G\to G$
is the identity map. In this situation, a version of the Fundamental
Theorem appears in \cite{Virelizier}, so let us compare this to our
\thref{3.3}. $X_G=\{\ol{e}\}$, the $G$-action on $G$ has only one orbit,
and
$$\ul{M}^{{\rm co}\ul{H}}_{\ol{e}}=\{(m_g)\in \prod_{g\in G}M_g~|~
\rho_{g,h}(m_{gh})=m_g\ot 1_h,~{\rm for~all~}g,h\in G\}$$
Virelizier then considers the image of the projection on the $g$-component
and calls this $M^{{\rm co}{\ul H}}_g$. Then he takes $M^{{\rm co}\ul{H}}$ to be
the product of all the $M^{{\rm co}{\ul H}}_g$. He then defines the tensor
product in a different way (see \cite[2.6]{Virelizier}), namely
$$(M^{{\rm co}\ul{H}}\ot \ul{H})_g= M^{{\rm co}{\ul H}}_g\ot H_g$$
while we have
$$(\ul{M}^{{\rm co}\ul{H}}\ot \ul{H})_g=\ul{M}^{{\rm co}\ul{H}}_{\ol{e}}\ot H_g$$
We will show next that these two tensor products are isomorphic, so that
Virelizier's Fundamental Theorem is a special case of \thref{3.3}
(and a fortiori of \cite[Theorem 3.4]{Takeuchi}).
\end{remark}

\begin{proposition}
Let $k$ be a field. With notation as in \reref{3.4}, we have that
$$M^{{\rm co}{\ul H}}_g\cong \ul{M}^{{\rm co}\ul{H}}_{\ol{e}},$$
in other words, the projection $p_g:\ \ul{M}^{{\rm co}\ul{H}}_{\ol{e}}
\to M^{{\rm co}{\ul H}}_g$ is injective.
\end{proposition}

\begin{proof}
Take $(m_h)_{h\in G}\in \ul{M}^{{\rm co}\ul{H}}_{\ol{e}}$. We have
to show that the full string is completely determined by one single entry
$m_g$. This can be seen as follows.
$$\rho_{h,h^{-1}g}(m_g)=m_h\ot 1_{h^{-1}g}\in M_h\ot H_{h^{-1}g}$$
Since $k$ is a field, the subspace of $H_{h^{-1}g}$ spanned by 
$1_{h^{-1}g}$ has a complement in $H_{h^{-1}g}$, so we have a projection $p$
$H_{h^{-1}g}\to k$, mapping $1_{h^{-1}g}$ to $1$. Applying $p$ to the second
tensor factor, we see that
$$m_h=(M_h\ot p)(\rho_{h,h^{-1}g}(m_g)).$$
\end{proof}

\section{Yetter-Drinfeld modules and the center construction}\selabel{4}
Let $H$ be a Hopf group-coalgebra. Yetter-Drinfeld modules over $H$ have been
studied in \cite{Zunino1,Zunino2,Wang1}. The aim of this Section is to recover
(a generalized version) of Yetter-Drinfeld modules using the centre construction.
This result is inspired by the classical result that, for a Hopf algebra $H$,
the centre of the category of $H$-modules is isomorphic to the category of
Yetter-Drinfeld modules. First we need some set theory. The following definition
goes back to Whitehead.

\begin{definition}\delabel{4.1}
Let $G$ be a monoid. A right crossed $G$-set is a right $G$-set $V$
together with a map $\nu:\ V\to G$ such that
$$g\nu(v\cdot g)=\nu(v)g,$$
for all $g\in G$ and $v\in V$. A morphism between two right crossed
$G$-sets $(V,\nu)$ and $(W,\omega)$ is a morphism of $G$-sets
$f:\ V\to W$ such that $\nu=\omega\circ f$.
The category of right crossed $G$-sets will be denoted by $\Xx_G^G$.
\end{definition}

The category $\Ss_G$ of right $G$-sets is a monoidal category. The tensor
product of two $G$-sets $X$ and $Y$ is the cartesian product, with the
diagonal action. The unit is $\{*\}$ with the trivial $G$-action
$*\cdot g=*$. The following result is classical.

\begin{proposition}\prlabel{4.2}
Let $G$ be a monoid.
The left weak center $\Ww_l(\Ss_G)$ is isomorphic to the category of right
crossed $G$-sets $\Xx_G^G$. If $G$ is a group, then it is also isomorphic to
the center.
\end{proposition}

\begin{proof} (sketch) 
Let $(V,s)\in \Ww_l(\Ss_G)$, and consider $s_G:\ V\times G\to G\times V$,
and write
$$s_G(v,1)=(\nu(v),f(v)),$$
with $\nu:\ V\to G$, $f:\ V\to V$.\\
Take a right $G$-set $X$, and consider the morphism of $G$-sets
$\varphi:\ G\to X$, $\varphi(g)=xg$. From the naturality of $s$, it follows
that
$$s_X\circ (V,\varphi)=(\varphi, V)\circ s_G.$$
Applying this to $(v,1)$, we see that
$s_X(v,x)=(x\nu(v),f(v))$. We then that
$$(*,v)=s_*(v,*)=(*\cdot \nu(v),f(v))=(*,f(v)),$$
so it follows that $f=V$ is the identity map on $V$. We conclude that
$s$ is completely determined by $\nu$:
\begin{equation}\eqlabel{4.2.1}
s_X(v,x)=(x\cdot \nu(v), v).
\end{equation}
Using the fact that $s_G$ is right $G$-linear, we find
$$(\nu(v)g,v\cdot g)=s_G(v,1)\cdot 1=s_G(v\cdot g,g)=(g\nu(v\cdot g),v\cdot g),$$
so $f(v)g=g\nu(v\cdot g)$, proving that $(V,\nu)$ is a right crossed $G$-set.\\
Conversely, if $(V,\nu)$ is a right crossed $G$-set, then we define
$(V,s)\in \Ww_l(\Ss_G)$ using \equref{4.2.1}.
\end{proof}

Let $G$ be a set. As we have seen, there is a unique coalgebra structure on
$G$ in the category of sets. We can then consider the category $\Ss^G$ of
right $G$-comodules. Its objects are couples $(X,f)$, with $X$ a set,
and $f:\ X\to G$ a map. A morphism between $(X,f)$ and $(Y,g)$ is a map
$\xi:\ X\to Y$ such that $g\circ \xi=f$. If $G$ is a monoid, then
$\Ss^G$ is a monoidal category. The tensor product is
$$(X,f)\times (Y,g)=(X\times Y, fg),$$
with $fg:\ X\times Y\to G$, $(fg)(x,y)=f(x)g(y)$. The unit is $(\{*\},i)$,
with $i(*)=1$.\\
For any set $X$, we can consider the right $G$-comodule $X_1=(M,1)$, with
$1:\ X\to G$ the constant map taking value $1$. Then for every 
$X=(X,f)\in \Ss^G$, the map
$$\rho:\ X\to X_1\times G,~~\rho(x)=(x,f(x)),$$
is a morphism in $\Ss^G$. This will be needed to prove the following result,
which is perhaps not so well-known as \prref{4.2}.

\begin{proposition}\prlabel{4.3}
Let $G$ be a monoid.
The left weak center $\Ww_l(\Ss^G)$ is isomorphic to the category of right
crossed $G$-sets $\Xx_G^G$. If $G$ is a group, then it is also isomorphic to
the center.
\end{proposition}

\begin{proof}
Take $((V,\nu),s)\in \Ww_l(\Ss^G)$. Take a set $X$, and fix $x\in X$. The map
$$g:\ \{*\}\to X_1,~~g(*)=x$$
is a morphism in $\Ss^G$. From the naturality of $s$, it follows that the
following diagram is commutative
$$\begin{diagram}
V\times\{*\}&\rTo^{s_*=V}&\{*\}\times V\\
\dTo^{(V,g)}&&\dTo_{(g,V)}\\
V\times X_1&\rTo^{s_{X_1}}&X_1\times V
\end{diagram}$$
and it follows that $s_{X_1}(v,x)=(x,v)$.\\
$G=(G,G)$ is an object of $\Ss^G$, hence we can consider
$s_G:\ V\times G\to G\times V$; let us denote
$$s_G(v,g)=(\gamma(v,g), v\cdot g),$$
for all $v\in G$ and $g\in G$. Since $((V,\nu),s)\in \Ww_l(\Ss^G)$, we have
that
$$(X_1,s_G)\circ (s_{X_1},G)=s_{X_1\times G},$$
implying that
$$s_{X_1\times G}(v,x,g)=(x,\gamma(v,g),v\cdot g).$$
From the naturality of $s$, we have the following commutative diagram
$$\begin{diagram}
V\times X&\rTo^{s_X}&X\times V\\
\dTo^{(V,\rho)}&&\dTo_{(\rho,V)}\\
V\times X_1\times G&\rTo^{s_{X_1\times G}}&X_1\times G\times V
\end{diagram}$$
Take $(v,x)\in V\times M$, and write $s_M(v,x)=(y,w)$. The commutativity
of the diagram implies that
$$(y,f(y),w)=s_{X_1\times G}(v,x,f(x))=(x,\gamma(v,f(x)),v\cdot f(x)).$$
Looking at the first component, we see that $y=x$. In particular, 
taking $M=G$, it follows that $\gamma(v,g)=g$. Looking at the third component,
we see that $w=v\cdot f(x)$, and we conclude that
\begin{equation}\eqlabel{4.3.1}
s_X(v,x)=(x,v\cdot f(x)),
\end{equation}
and, in particular,
$$s_G(v,g)=(g,v\cdot g).$$
$s_G$ is a morphism in $\Ss^G$, hence $G\nu\circ s_G=\nu G$, and it follows that
$$\nu(v)g=g \nu(v\cdot g),$$
for all $g\in G$ and $v\in V$. From the fact that
$$s_{G\times G}=(G,s_V)\circ (s_V,G),$$
we conclude that $s_{G\times G}(v,(g,h))=((g,h), (v\cdot g)\cdot h)$. Since
$m:\ G\times G\to G$ is a morphism in $\Ss^G$, we have the following commutative
diagram from the naturality of $s$:
$$\begin{diagram}
V\times G\times G&\rTo^{s_{G\times G}}&G\times G\times V\\
\dTo^{(V,m)}&&\dTo_{(m,V)}\\
V\times G&\rTo^{s_{G}}&G\times V
\end{diagram}$$
and it follows that $(v\cdot g)\cdot h=v\cdot (gh)$. Finally, the constant
map $1:\ \{*\}\to G$ is a morphism in $\Ss^G$, so we have a commutative
diagram
$$\begin{diagram}
V\times \{*\}&\rTo^{s_{*}}&\{*\}\times V\\
\dTo^{(V,1)}&&\dTo_{(1,V)}\\
V\times G&\rTo^{s_{G}}&G\times V
\end{diagram}$$
from which we deduce that $v\cdot 1=v$. Thus we have shown that $(V,\nu)$ is a
right crossed $G$-module.\\
Conversely, given a right crossed $G$-module $(V,\nu)$, we define
$((V,\nu),s)\in \Ww_l(\Ss^G)$ using \equref{4.3.1}.
\end{proof}

Now let $\ul{H}=(G,(H_g)_{g\in G})$ be a semi-Hopf group coalgebra. In particular,
this means that $G$ is a monoid.

\begin{definition}\delabel{4.4}
A right-right $\ul{H}$-Yetter-Drinfeld module consists of the following data:
\begin{itemize}
\item An object $\ul{M}=(V,(M_v)_{v\in V})\in \Tt_k$;
\item a right $\ul{H}$-comodule structure $\ul{\rho}:\ \ul{M}\to \ul{M}\ot
\ul{H}$; in particular, $V$ is a right $G$-set;
\item a right $\ul{H}$-module structure $\ul{\psi}:\ \ul{M}\ot
\ul{H}\to \ul{M}$; this implies that we have a map $\nu:\ V\to V\times G$
making $(V,f)\in \Ss^G$,
\end{itemize}
satisfying the following compatibility relations:
\begin{itemize}
\item $V$ is a crossed right $G$-set;
\item for all $v\in V$, $g\in G$, $m\in M_{v\cdot g}$ and
$h\in H_{g\nu(v\cdot g)}=H_{\nu(v)g}$, we have
\begin{equation}\eqlabel{4.4.1}
\bigl(mh_{(2,\nu(v\cdot g))}\bigr)_{[0,v]}\ot
h_{(1,g)}\bigl(mh_{(2,\nu(v\cdot g))}\bigr)_{[1,g]}=
m_{[0,v]}h_{(1,\nu(v))}\ot m_{[1,g]}h_{(2,g)}.
\end{equation}
\end{itemize}
The category of right-right Yetter-Drinfeld modules and $\ul{H}$-linear
$\ul{H}$-colinear maps is denoted by $\YD_{\ul{H}}^{\ul{H}}$.
\end{definition}

The category of right $\ul{H}$-modules is a monoidal category.
For $\ul{P}=(X,P),\ul{Q}=(Y,Q)\in \Tt_{\ul{H}}$, we have the following
right $\ul{H}$-coaction on $\ul{P}\ot \ul{Q}$:
$$\begin{matrix}
\ul{P}\ot\ul{Q}\ot\ul{H}&\mapright{\ul{\psi}}&\ul{P}\ot\ul{Q}\\
X\times Y\times G&\mapleft{X\times Y\times fg}& X\times Y\\
P_x\ot Q_y\ot H_{f(x)g(y)}&\mapright{\psi_{x,y}}& P_x\ot Q_y,
\end{matrix}$$
with
$$\psi_{x,y}(p\ot q\ot h)=ph_{(1,f(x))}\ot qh_{(1,g(y))},$$
for all $p\in P_x$, $q\in Q_x$ and $h\in H_{f(x)g(y)}$. We can therefore
consider the centre and the weak centre of $\Tt_{\ul{H}}$, and 
the main result of this Section is now the following.

\begin{theorem}\thlabel{4.5}
Let $\ul{H}$ be a semi-Hopf group coalgebra. Then we have an isomorphism
of categories
$$\Ww_r(\Tt_{\ul H})\cong \YD_{\ul{H}}^{\ul{H}}.$$
If $\ul{H}$ is a Hopf group coalgebra, then the weak centre is equal to
the centre.
\end{theorem}

\begin{proof}
Take $(\ul{M},\ul{t})\in \Ww_r(\Tt_{\ul{H}})$. For every 
$\ul{P}=(X,P)\in \Tt_{\ul{H}}$, we have a morphism
$$\begin{matrix}
\ul{P}\ot \ul{M}&\mapright{\ul{t}_{\ul{P}}}&\ul{M}\ot \ul{P}\\
X\times V&\mapleft{s_X}&V\times X\\
P_x\ot M_{v\cdot f(x)}&\mapright{t_{P,v,x}}&M_v\ot P_x,
\end{matrix}$$
with $s_X(v,x)=(x,v\cdot f(x))$. Since the forgetful functor
$\Tt_k\to \dul{\rm Sets}^{\rm op}$ is strongly monoidal, it follows that
$(V, s)\in \Ww_l(\Ss^G)$, so $V$ is a right crossed $G$-module, by
\prref{4.3}. Look at $t_{\ul{H}}$, and consider the composition
$$\ul{\rho}=t_{\ul{H}}\circ (\ul{\eta}\ot \ul{M}),$$
namely
$$\begin{matrix}
\ul{M}&\mapright{\ul{\rho}}&\ul{M}\ot \ul{H}\\
V&\mapleft{}&V\times G\\
M_{v\cdot g}&\mapright{\rho_{v,g}}&M_v\ot H_g
\end{matrix}$$
We will use the following Sweedler-type notation:
$$\rho_{v,g}(m)=m_{[0,v]}\ot m_{[1,g]}.$$
\ul{step 1}. $t$ is completely determined by $\rho$.\\
For each $x\in X$, we fix $p_x\in P_x$, and consider the following
morphism in $\Tt_{\ul H}$:
$$\begin{matrix}
\ul{H}&\mapright{\ul{\varphi}}&\ul{P}\\
G&\mapleft{f}& X\\
H_{f(x)}&\mapright{\varphi_x}&P_x,
\end{matrix}$$
with $\varphi_x(h)=p_x h$, for all $h\in H_{f(x)}$.
From the naturality of $\ul{t}$, it follows that the following diagram
is commutative:
$$\begin{diagram}
\ul{H}\ot \ul{M}&\mapright{\ul{t}_{\ul{H}}}& \ul{M}\ot \ul{H}\\
\dTo^{\ul{f}\ot\ul{V}}&&\dTo_{\ul{V}\ot\ul{f}}\\
\ul{P}\ot \ul{M}&\mapright{\ul{t}_{\ul{P}}}& \ul{M}\ot \ul{P}.
\end{diagram}$$
In particular, we have a commutative diagram, for all $v\in V$ and
$x\in X$:
$$\begin{diagram}
H_{f(x)}\ot M_{v\cdot f(x)}&\mapright{{t}_{\ul{H},v,f(x)}}&
M_v\ot H_{f(x)}\\
\dTo^{\varphi_x\ot M_{v\cdot f(x)}}&&\dTo_{M_v\ot \varphi_x}\\
P_x\ot M_{v\cdot f(x)}&\mapright{{t}_{\ul{P},v,x}}&
M_v\ot P_x
\end{diagram}$$
and it follows that
\begin{equation}\eqlabel{4.5.1}
t_{\ul{P},v,x}(p_x \ot m)=m_{[0,v]}\ot p_x m_{[1,f(x)]}.
\end{equation}
\ul{step 2}. From the definition of the right center, it follows that the
following diagram commutes:
$$\begin{diagram}
\ul{H}\ot\ul{H}\ot \ul{M}&\mapright{\ul{t}_{\ul{H}\ot\ul{H}}}& \ul{M}\ot \ul{H}\ot\ul{H}\\
\SE_{\ul{H}\ot\ul{t}_{\ul{H}}}&&\NE_{\ul{t}_{\ul{H}}\ot \ul{H}}\\
&\ul{H}\ot\ul{M}\ot \ul{H}&
\end{diagram}$$
It follows that we have the following commutative diagram, for all
$g,h\in G$ and $v\in V$:
$$\begin{diagram}
H_g\ot H_h\ot M_{v\cdot(gh)}&\mapright{{t}_{\ul{H}\ot\ul{H},v,gh}}&
M_v\ot H_g\ot H_h\\
\SE_{H_g\ot {t}_{\ul{H},v\cdot g, h}}&&\NE_{{t}_{\ul{H},v,g}\ot H_h}\\
&H_g\ot M_{v\cdot g}\ot H_h
\end{diagram}$$
Evaluating this diagram at $1_g\ot 1_h\ot m$, with $m\in M_{v\cdot(gh)}$,
we find, using \equref{4.5.1}
$$m_{[0,v]}\ot \Delta_{g,h}(m_{[1,gh]})=
\rho_{v,g}(m_{[0,v]})\ot m_{[1,gh]},$$
so $\ul{\rho}$ is coassociative.\\
\ul{step 3}. It follows from the definition of the left center that
$\ul{t}_{\ul{k}}=\ul{M}:\ \ul{k}\ot \ul{M}\to \ul{M}\ot \ul{k}$. Using
\equref{4.5.1}, we compute, for all $m\in M_v$ that
$$m=t_{\ul{k},v,*}(m)=m_{[0,v]}\ot 1\cdot m_{[1,1]}=m_{[0,v]}\varepsilon(m_{[1,1]}),$$
and the counit property of $\ul{\rho}$ follows. Thus $\ul{\rho}$
defines a right $\ul{H}$-comodule structure on $\ul{M}$.\\
\ul{step 3}. $\ul{t}_{\ul{H}}$ is right $\ul{H}$-linear, hence the
following diagram commutes:
$$\begin{diagram}
\ul{H}\ot\ul{M}\ot\ul{H}&\mapright{}&\ul{H}\ot\ul{M}\\
\dTo^{\ul{t}_{\ul{H}}\ot \ul{H}}&&\dTo_{\ul{t}_{\ul{H}}}\\
\ul{M}\ot\ul{H}\ot\ul{H}&\mapright{}&\ul{M}\ot\ul{H}
\end{diagram}$$
For all $v\in V$ and $g\in G$, we have the following commutative diagram,
keeping in mind that $g\nu(v\cdot g)=\nu(v)g$:
$$\begin{diagram}
H_g\ot M_{v\cdot g}\ot H_{g\nu(v\cdot g)}&\mapright{}&H_g\ot M_{v\cdot g}\\
\dTo^{t_{\ol{H},v,g}\ot H_{g\nu(v\cdot g)}}&&\dTo_{t_{\ol{H},v,g}}\\
M_v\ot H_g \ot H_{\nu(g)g}&\mapright{}&M_v\ot H_g
\end{diagram}$$
Evaluating this diagram at $1_g\ot m\ot h$, with $m\in M_{v,g}$,
$h\in H_{g\nu(v\cdot g))}$, we find that \equref{4.4.1} holds. Hence
it follows that $\ul{M}$ is a right-right Yetter-Drinfeld module.\\

Conversely, let $\ul{M}$ be a right-right Yetter-Drinfeld module. We define
$$\ul{t}:\ -\ot \ul{M}\to \ul{M}\ot -$$
using \equref{4.5.1}. Standard computations show that
$(\ul{M},\ul{t})\in \Ww_l(\Tt_{\ul{H}})$.
\end{proof}

In a similar way, we can compute the right weak center of ${}_{\ul{H}}\Tt$,
and this is isomorphic to the category of left-right Yetter-Drinfeld
modules ${}_{\ul{H}}\YD^{\ul{H}}$. Objects in ${}_{\ul{H}}\YD^{\ul{H}}$
are Turaev $k$-modules $\ul{M}$ with a right $\ul{H}$-comodule structure
$\ul{\rho}:\ \ul{M}\to \ul{M}\ot \ul{H}$ and a left $\ul{H}$-module structure
$\ul{\psi}:\ \ul{H}\ot \ul{M}\to \ul{M}$ such that
$$(\ul{\psi}\ot\ul{\mu})\circ (\ul{H}\ot\ul{\tau}\ot \ul{H})\circ
(\ul{\Delta}\ot \ul{\rho})=
(\ul{M}\ot\ul{\mu})\circ(\ul{\rho}\ot\ul{H})\circ\ul{\tau}\circ
(\ul{H}\ot\ul{\psi})\circ(\ul{\Delta}\ot\ul{M}).$$

\begin{example}\exlabel{4.6}
Let $G$ be a group, fix $g\in G$, and let $V$ be the orbit of $g$
under the adjoint action $h\triangleleft k=k^{-1}hk$:
$$V=\{k^{-1}gk~|~k\in G\}.$$
Let $\nu:\ V\ot G$ be the embedding of $V$ into $G$. Then $(V,\nu)$, with
the adjoint action, is a right crossed $G$-set. Assume that $\ul{H}$ is a
crossed semi-Hopf $G$-coalgebra. Recall from \cite{Turaev,Zunino1,Zunino2}
that this means that $\ul{H}$ is a semi-Hopf $G$-coalgebra with a family
of algebra isomorphisms
$$\varphi^k_h:\ H_k\to H_{hkh^{-1}},$$
satisfying the following conditions (we omit the upper index if no confusion is
possible):
$$\varphi_h\circ \varphi_k=\varphi_{kh}~~;~~\varphi^k_1=H_k;$$
$$(\varphi_k\ot \varphi_k)\circ\Delta_{l,h}=\Delta_{klk^{-1},khk^{-1}}~~;~~
\varepsilon\circ \varphi^1_g=\varepsilon.$$
We also assume that the following additional condition holds:
\begin{equation}\eqlabel{4.6.1}
\varphi^l_h=\varphi^l_k~~{\rm if}~~hlh^{-1}=klk^{-1}.
\end{equation}
Take $M\in {}_{H_g}\Mm$. For each $v=h^{-1}gh\in V$, let $M_v=M$ as
a $k$-module, with left $H_v$-action $a\cdot m=\varphi^v_h(a)m$, for all
$a\in H$. This is well-defined, because of \equref{4.6.1}, and
$\ul{M}=(V,(M_v)_{v\in V})\in {}_{\ul{H}}\Tt$.\\
Assume now that we have a family of maps $\rho_l:\ M\to M\ot H_l$, indexed
by $l\in G$, as in \exref{3.2.0}. Then we have the following right $\ul{H}$-coaction
on $\ul{M}$:
$$\rho_{v,k}=\rho_k:\ M_{k^{-1}vk}=M\to M_v\ot H_k=M\ot H_k.$$
We use the Sweedler notation
$$\rho_l(m)=m_{[0]}\ot m_{[1,l]}.$$
With this action and coaction, $\ul{M}\in {}_{\ul{H}}\YD^{\ul{H}}$ if and only if
\begin{eqnarray*}
&&\hspace*{-2cm}
\varphi_h(a_{(1,h^{-1}gh}))m_{[0]}\ot a_{(2,l)}m_{[1,l]}\\
&=&
\bigl(\varphi_{hl}(a_{(2,l^{-1}h^{-1}ghl)})m\bigr)_{[0]}\ot
\bigl(\varphi_{hl}(a_{(2,l^{-1}h^{-1}ghl)})m\bigr)_{[1,l]}a_{(1,l)},
\end{eqnarray*}
for all $h,l\in G$ and $a\in H_{h^{-1}ghl}$. Taking $h=1$, we find
\begin{eqnarray}
&&\hspace*{-2cm}
a_{(1,g}))m_{[0]}\ot a_{(2,l)}m_{[1,l]}\nonumber\\
&=&
\bigl(\varphi_{l}(a_{(2,l^{-1}gl)})m\bigr)_{[0]}\ot
\bigl(\varphi_{l}(a_{(2,l^{-1}gl)})m\bigr)_{[1,l]}a_{(1,l)},\eqlabel{4.6.2}
\end{eqnarray}
for all $l\in G$ and $a\in H_{gl}$.\\
Recall the definition of $g$-Yetter-Drinfeld module, introduced by Zunino
in \cite{Zunino2}. Let $\ul{H}$ be a crossed Hopf $G$-coalgebra, and fix $g\in G$
and $M\in {}_{H_g}\Mm$. Assume that we have maps  $\rho_l$ as above.
Then $M$ is a $g$-Yetter-Drinfeld module if the following compatibility conditions
are satisfied, for all $l\in G$ and $a\in H_{gl}$:
\begin{eqnarray}
&&\hspace*{-2cm}
a_{(1,g}))m_{[0]}\ot a_{(2,l)}m_{[1,l]}\nonumber\\
&=&
(a_{(2,g)}m)_{[0]}\ot (a_{(2,g)}m)_{[1,l]}\varphi_{\alpha^{-1}}(a_{(1,glg^{-1})}).\eqlabel{4.6.3}
\end{eqnarray}
Observe that the conditions \equref{4.6.2} and \equref{4.6.3} are not the same.
\end{example}

\section{Further generalizations}\selabel{5}
\begin{blanco} In the definition of the Turaev and Zunino categories, we
can replace the category of $k$-modules by any monoidal category $\Cc$.
The Turaev category $\Tt_{\Cc}$ and the Zunino category $\Zz_{\Cc}$ will
be symmetric (resp. braided) if $\Cc$ is symmetric (resp. braided).
\end{blanco}

\begin{blanco} The category $\dul{\rm Sets}$ can be replaced by any
full subcategory of $\dul{\rm Sets}$ containing a singleton and closed
under finite cartesian products, for example the category of finite sets,
or the category of countable sets.
\end{blanco}

\begin{blanco} (E. Villanueva)
Take $\ul{M}=(X,(M_x)_{x\in X})\in \Tt_k$, and consider $X$ as a topological
space with the discrete topology. $\ul{M}$ can then be viewed as a
(pre)sheaf of $k$-modules on $X$. Let $\Ss_k$ be the category with
objects $(X,\Ff)$, with $X$ a topological space, and $\Ff$ a sheaf of
$k$-modules on $X$. Then $\Tt_k$ is a full subcategory of $\Ss_k$.
\end{blanco}

\end{document}